\documentclass[12pt]{amsart}
\usepackage{amsmath,amssymb,amsbsy,amsfonts,amsthm,latexsym,amsopn,amstext,amsxtra,euscript,
amscd}
\usepackage{hyperref}
\usepackage{color}
\begin{document}

\title[Sieve Procedure for the M\"obius prime-functions ...]{Sieve Procedure\\
for the M\"obius prime-functions,\\
the Infinitude of Primes\\
and the Prime Number Theorem}

\author[R. M.~Abrarov]{R. M.~Abrarov$^1$}
\email{$^1$rabrarov@gmail.com}
\author[S. M.~Abrarov]{S. M.~Abrarov$^2$}
\email{$^2$absanj@gmail.com}

\date{September 27, 2011}

\begin{abstract}

Using a sieve procedure akin to the sieve of Eratosthenes we show how for each prime $p$ to build the corresponding M\"obius prime-function, which in the limit of infinitely large primes becomes identical to the original M\"obius function. Discussing this limit we present two simple proofs of the Prime Number Theorem. In the framework of this approach we give several proofs of the infinitude of primes.

\end{abstract}

\maketitle

\section{\textbf{Sieve procedure}}

\textbf{Some definitions}. A number is squarefree \cite{Weisstein3} if in the~prime factorization of~this number no~prime factor occurs more than once. The M\"obius function \textit{$\mu \left(n\right)$} \cite{Hardy, Weisstein1} is defined for positive integers $n$ by
\begin{equation} \label{Eq_1}
\begin{cases} 
\mu \left(1\right)=1,\\
\mu \left(n\right)=1,  &\text{\scriptsize if $n$ is squarefree with an even number of prime factors,}\\
\mu \left(n\right)=-1, &\text{\scriptsize if $n$ is squarefree with an odd number of prime factors,}\\
\mu \left(n\right)=0,  &\text{\scriptsize if $n$ is not squarefree.}
\end{cases}
\end{equation}
A number is squareful \cite{Weisstein3} (opposite, complementary for squarefree) if its prime decomposition contains at least one square (if and only if \textit{$\mu \left(n\right)=0$}). Direct result of these definitions is that all prime numbers $p$ are squarefree with \textit{$\mu \left(p\right)=-1$}.
 
The Mertens function $M\left(n\right)$ \cite{Hardy, Wikipedia1} is the cumulative sum of the M\"obius function
\begin{equation} \label{Eq_2} 
M\left(n\right)=\sum _{k=1}^{n}\mu \left(k\right).         
\end{equation}

There is an identity \cite{RAbrarov1, WolframDemo} between the M\"obius function and the delta functions
\begin{equation} \label{Eq_3} 
\mu \left(n\right)=-\sum _{i,j=1}^{\sqrt{n} }\mu \left(i\right)\mu \left(j\right)\,  \delta \left(\frac{n}{i\cdot j} \right), \qquad    n\ge 2,     \end{equation} 
where the delta function
$$\delta\left(x\right)=
\begin{cases} 
1, \text{if $x \in \mathbb{N}_0=\{0, 1, 2,\dots\}$}  \\
0, \text{if $x \notin \mathbb{N}_0$}
\end{cases} .$$

Define the prime detecting function with outcomes 1 for primes and 0 for composite numbers as
\footnotesize
\begin{equation} \label{Eq_4} 
\begin{aligned}
p^{\delta} \left(n\right) & = {\left(1-\delta \left(\frac{n-2^{2} }{2} \right)\right)\left(1-\delta \left(\frac{n-3^{2} }{3} \right)\right)\left(1-\delta \left(\frac{n-5^{2} }{5} \right)\right)...} \\ 
& = {\prod_{\text{$p$ prime}}\left(1-\delta \left(\frac{n-p^{2} }{p} \right)\right),\qquad n\ge 2\; .} \end{aligned}
\end{equation}
\normalsize
If we assign $\pi \left(n\right)$ the number of primes $p\le n$ then
\begin{equation} \label{Eq_5} 
\pi \left(n\right)=\sum _{m=2}^{n}p^{\delta} \left(m\right) .
\end{equation}

Let us denote $Q\left(n\right)$, $Q^{+} \left(n\right)$, and $Q^{-} \left(n\right)$ accordingly as the number of squarefree integers, the number of integers with positive values, and the number of integers with negative values of the M\"obius function between 1 and~\textit{$n$}. From these definitions we have
\begin{equation} \label{Eq_6} 
\begin{array}{rcl} 
{Q\left(n\right)} & {=} & {Q^{+} \left(n\right)+Q^{-} \left(n\right),} \\ 
{M\left(n\right)} & {=} & {Q^{+} \left(n\right)-Q^{-} \left(n\right).} 
\end{array} 
\end{equation} 
The corresponding asymptotic densities are
$$
\mathop{\lim }\limits_{n\to \infty } \frac{Q\left(n\right)}{n} , \mathop{\lim }\limits_{n\to \infty } \frac{M\left(n\right)}{n} =\mathop{\lim }\limits_{n\to \infty } \left(\frac{Q^{+} \left(n\right)}{n} -\frac{Q^{-} \left(n\right)}{n} \right).
$$
We can define also the asymptotic densities above differently.  If we randomly choose $n$ integers between 1 and \textit{$n$}, then the probability $\Pr$ that number is squarefree, or the integer with positive or negative values of the M\"obius function can be determined by the frequencies
$$
\frac{Q\left(n\right)}{n} , \frac{Q^{+} \left(n\right)}{n} , \frac{Q^{-} \left(n\right)}{n} ,
$$ 
where $Q\left(n\right)$, $Q^{+} \left(n\right)$, and $Q^{-} \left(n\right)$ are the quantity of events with nonzero, positive and negative values of the M\"obius function, respectively. If these frequencies tend to some limits as $n$ tends to infinity, we call these limits by corresponding asymptotic densities.

Introduce $\tilde{\mu }_{p} \left(n\right)$, where $p$ either 1 or prime - a set of the modified M\"obius functions - and call them the M\"obius prime-functions or simply the M\"obius $p$-functions. Consequently, the corresponding modified functions are
\begin{equation} \label{Eq_7} 
\tilde{M}_{p} \left(n\right)=\sum _{k=1}^{n}\tilde{\mu }_{p} \left(k\right) =\tilde{Q}_{p}^{+} \left(n\right)-\tilde{Q}_{p}^{-} \left(n\right) 
\end{equation} 
- the Mertens $p$-functions, and
\begin{equation} \label{Eq_8} 
\tilde{Q}_{p} \left(n\right)=\sum _{k=1}^{n}\left|\tilde{\mu }_{p} \left(k\right)\right| =\tilde{Q}_{p}^{+} \left(n\right)+\tilde{Q}_{p}^{-} \left(n\right)   .       
\end{equation} 
$\tilde{Q}_{p} \left(n\right)$ is ''squarefree'' counting $p$-function (we use quotes because the numbers are squarefree relative to the M\"obius $p$-function $\tilde{\mu }_{p} \left(n\right)$; the actual squarefree number is defined in \eqref{Eq_1} by the original M\"obius function $\mu \left(n\right)$).

Define first $\tilde{\mu }_{1} \left(n\right)$ as
\begin{equation} \label{Eq_9} 
\begin{array}{l} 
{\tilde{\mu }_{1} \left(1\right)=1,} \\ 
{\tilde{\mu }_{1} \left(n\right)=-\delta \left(n-1^{2} \right)=-1,\qquad n\ge 2,} 
\end{array} 
\end{equation} 
i.e. in this definition in the first approximation all natural numbers except 1 we treat as primes. Below are the values of $\tilde{\mu }_{1} \left(n\right)$ for the first 30 natural numbers:
$$
\setcounter{MaxMatrixCols}{30}
\begin{smallmatrix}
1&2&3&4&5&6&7&8&9&10&11&12&13&14&15&16&17&18&19&20&21&22&23&24&25&26&27&28&29&30\\
+&-&-&-&-&-&-&-&-&-&-&-&-&-&-&-&-&-&-&-&-&-&-&-&-&-&-&-&-&-
\end{smallmatrix}\dots
$$
We will refer to this infinite row as the Complete Sequence. Subsequently, we will make corrections in the Complete Sequence as soon as we proceed to the next $\tilde{\mu }_{p} \left(n\right)$.
The difference between the asymptotic densities of positive and negative values of  $\tilde{\mu }_{1} \left(n\right)$ is
\footnotesize
\begin{equation} \label{Eq_10} 
\mathop{\lim }\limits_{n\to \infty } \frac{\tilde{M}_{1} \left(n\right)}{n} =\mathop{\lim }\limits_{n\to \infty } \left(\Pr \left(\tilde{\mu }_{1} \left(n\right)=+1\right)-\Pr \left(\tilde{\mu }_{1} \left(n\right)=-1\right)\right)=-1.    
\end{equation} 
\normalsize
Corresponding ''squarefree'' asymptotic density for $\tilde{Q}_{1} \left(n\right)$ is
\begin{equation} \label{Eq_11} 
\mathop{\lim }\limits_{n\to \infty } \frac{\tilde{Q}_{1} \left(n\right)}{n} =\mathop{\lim }\limits_{n\to \infty } \left(\Pr \left(\left|\tilde{\mu }_{1} \left(n\right)\right|=1\right)\right)=1.
\end{equation} 
To get the next M\"obius $p$-function we pick the number following 1 in the Complete Sequence. It is prime 2 and the corresponding M\"obius $p$-function is $\tilde{\mu }_{2} \left(n\right)$. For $n\ge 2$ define
\begin{equation*}
\begin{array}{l}
\tilde{\mu }_{2} \left(2\cdot n\right)=+1 \qquad \text{if 2$\nmid n$},\\ 
\tilde{\mu }_{2} \left(2\cdot n\right)=0 \qquad \;\;\, \text{if 2$\mid n$}.
\end{array}
\end{equation*}
For all other $n$ $\tilde{\mu }_{2} \left(n\right)$ remains unchanged
\[\tilde{\mu }_{2} \left(n\right)=\tilde{\mu }_{1} \left(n\right).\] 
The altered Complete Sequence with values of the new function $\tilde{\mu }_{2} \left(n\right)$ for the first 30 numbers is
$$
\setcounter{MaxMatrixCols}{30}
\begin{smallmatrix}
1&2&3&4&5&6&7&8&9&10&11&12&13&14&15&16&17&18&19&20&21&22&23&24&25&26&27&28&29&30\\
+&-&-&0&-&+&-&0&-& +& -& 0&- &+ &- &0 &- &+ &- &0 &- &+ &- &0 &- &+ &- &0 &- &+
\end{smallmatrix}\dots
$$
Obviously, $\tilde{\mu }_{2} \left(n\right)$ can be expressed as
\begin{equation} \label{Eq_12}
\tilde{\mu }_{2} \left(n\right)=\tilde{\mu }_{1} \left(n\right)+2\delta \left(\frac{n-2^{2} }{2} \right)-\delta \left(\frac{n-2^{2} }{2^{2} } \right).
\end{equation}

Replacing the delta functions in \eqref{Eq_12} by corresponding frequencies or probabilities \cite{RAbrarov2}, we get the asymptotic densities
\begin{equation} \label{Eq_13}
\begin{array}{rcl} 
{\mathop{\lim }\limits_{n\to \infty } \frac{\tilde{M}_{2} \left(n\right)}{n} } & {=} & {\mathop{\lim }\limits_{n\to \infty } \left(\Pr \left(\tilde{\mu }_{2} \left(n\right)=+1\right)-\Pr \left(\tilde{\mu }_{2} \left(n\right)=-1\right)\right)} \\ {} & {=} & {-1+2\frac{1}{2} -\frac{1}{2^{2} } =-\left(1-\frac{1}{2} \right)^{2} }
\end{array}
\end{equation} 
and 
\begin{equation} \label{Eq_14} 
\mathop{\lim }\limits_{n\to \infty } \frac{\tilde{Q}_{2} \left(n\right)}{n} =\mathop{\lim }\limits_{n\to \infty } \left(\Pr \left(\left|\tilde{\mu }_{2} \left(n\right)\right|=1\right)\right)=1-\frac{1}{2^{2} }\; .
\end{equation}

In order to determine the next prime we sieve all odd numbers (in other words, we remove all multiples of 2)
$$
\setcounter{MaxMatrixCols}{30}
\begin{smallmatrix}
1&\quad\!\!&3&\quad\!\!&5&\quad\!\!&7&\quad\!\!&9&\quad\!\!&11&\quad\!\!&13&\quad\!\!&15&\quad\!\!&17&\quad\!\!&19&\quad\!\!&21&\quad\!\!&23&\quad\!\!&25&\quad\!\!&27&\quad\!\!&29&\\
+&&-&&-&&-&&-&& -&&- &&- &&- &&- &&- &&- &&- &&- &&- &
\end{smallmatrix}\dots
$$
In this new sequence we again pick the next number following 1. Now it is prime 3. Continue to the next M\"obius $p$-function $\tilde{\mu }_{3} \left(n\right)$ (\textit{n }runs only through the sifted numbers $(n,2)=1$ except 1) 
\begin{equation*}
\begin{array}{l}
\tilde{\mu }_{3} \left(3\cdot n\right)=+1 \qquad \text{if 3$\nmid n$},\\ 
\tilde{\mu }_{3} \left(3\cdot n\right)=0 \qquad \;\;\, \text{if 3$\mid n$}.
\end{array}
\end{equation*}
Now the new values in the sequence above for prime 3 (3-sequence) are
$$
\setcounter{MaxMatrixCols}{30}
\begin{smallmatrix}
1&\quad\!\!&3&\quad\!\!&5&\quad\!\!&7&\quad\!\!&9&\quad\!\!&11&\quad\!\!&13&\quad\!\!&15&\quad\!\!&17&\quad\!\!&19&\quad\!\!&21&\quad\!\!&23&\quad\!\!&25&\quad\!\!&27&\quad\!\!&29&\\
+&&-&&-&&-&&0&& -&&- &&+ &&- &&- &&+ &&- &&- &&0 &&- &
\end{smallmatrix}\dots
$$
Insert these new values into the Complete Sequence making at the same time the following changes ($n$ runs only through the sifted numbers $(n,2)=1$ except 1)
\begin{equation*}
\begin{array}{l}
\tilde{\mu }_{3} \left(2\cdot 3\cdot n\right)=-1 \qquad \text{if 3$\nmid n$},\\ 
\tilde{\mu }_{3} \left(2\cdot 3\cdot n\right)=0 \qquad \;\;\, \text{if 3$\mid n$}.
\end{array}
\end{equation*}
For all other $n$ $\tilde{\mu }_{3} \left(n\right)$ remains unchanged
\[\tilde{\mu }_{3} \left(n\right)=\tilde{\mu }_{2} \left(n\right)   .\] 
After these alterations the Complete Sequence becomes
$$
\setcounter{MaxMatrixCols}{30}
\begin{smallmatrix}
1&2&3&4&5&6&7&8&9&10&11&12&13&14&15&16&17&18&19&20&21&22&23&24&25&26&27&28&29&30\\
+&-&-&0&-&+&-&0&0& +& -& 0&- &+ &+ &0 &- &0 &- &0 &+ &+ &- &0 &- &+ &0 &0 &- &-
\end{smallmatrix}\dots
$$

All changes made above can be reflected as
\begin{equation} \label{Eq_15} 
\begin{array}{rcl} {\tilde{\mu }}_{3} {\left(n\right)} & {=} & {\tilde{\mu }_{2} \left(n\right)+2\underbrace{\left(1-\delta \left(\frac{n}{2} \right)\right)}_{{\rm sieving\; factor}} \delta \left(\frac{n-3^{2} }{3} \right)} 
\\ {} & {} & -2\underbrace{\left(1-\delta \left(\frac{n-2^{2} \cdot 3^{2} }{2^{2} } \right)\right)}_{{\rm sieving\; factor}} \delta \left(\frac{n-2\cdot 3^{2} }{2\cdot 3} \right)
\\ {} & {} & -\left(1-\delta \left(\frac{n}{2} \right)\right)\delta \left(\frac{n-3^{2} }{3^{2} } \right)
\\ {} & {} & +\left(1-\delta \left(\frac{n-2^{2} \cdot 3^{2} }{2^{2} } \right)\right)\delta \left(\frac{n-2\cdot 3^{2} }{2\cdot 3^{2} } \right). 
\end{array} 
\end{equation}
Here and further for simplicity and clarity we use the products of the delta functions. To return to the expression with single delta functions one can apply the plain rules
\begin{equation} \label{Eq_16} 
\begin{array}{l} {\delta \left(\frac{n}{r} \right)\delta \left(\frac{n-s^{2} }{s} \right)=\delta \left(\frac{n-r\cdot s^{2} }{r\cdot s} \right),} \\ 
{\delta \left(\frac{n-r^{2} \cdot s^{2} }{r^{2} } \right)\delta \left(\frac{n-s^{2} }{r\cdot s} \right)=\delta \left(\frac{n-r^{2} \cdot s^{2} }{r^{2} \cdot s} \right),} 
\end{array} 
\end{equation} 
where $r{\rm \; and\; }s$ are relatively prime (coprime) numbers.

Again replacing the delta functions in \eqref{Eq_15} by corresponding frequencies or probabilities, we get the asymptotic densities
\begin{equation} \label{Eq_17} 
\begin{array}{rcl} {\mathop{\lim }\limits_{n\to \infty } \frac{\tilde{M}_{3} \left(n\right)}{n} } & {=} & {\mathop{\lim }\limits_{n\to \infty } \left(\Pr \left(\tilde{\mu }_{3} \left(n\right)=+1\right)-\Pr \left(\tilde{\mu }_{3} \left(n\right)=-1\right)\right)} \\ 
{} & {=} & {-\left(1-\frac{1}{2} \right)^{2} \left(1-\frac{1}{3} \right)^{2} },
\end{array}    
\end{equation} 
and
\begin{equation} \label{Eq_18} 
\mathop{\lim }\limits_{n\to \infty } \frac{\tilde{Q}_{3} \left(n\right)}{n} =\mathop{\lim }\limits_{n\to \infty } \left(\Pr \left(\left|\tilde{\mu }_{3} \left(n\right)\right|=1\right)\right)=\left(1-\frac{1}{2^{2} } \right)\left(1-\frac{1}{3^{2} } \right).    
\end{equation}

Removing from the last sieved sequence (3-sequence) multiples of 3
$$
\setcounter{MaxMatrixCols}{30}
\begin{smallmatrix}
1&\quad\!\!&\quad\!\!&\quad\!\!&5&\quad\!\!&7&\quad\!\!&\quad\!\!&\quad\!\!&11&\quad\!\!&13&\quad\!\!&\quad\!\!&\quad\!\!&17&\quad\!\!&19&\quad\!\!&\quad\!\!&\quad\!\!&23&\quad\!\!&25&\quad\!\!&\quad\!\!&\quad\!\!&29&\\
+&&&&-&&-&&&& -&&- &&&&- &&- &&&&- &&- &&&&- &
\end{smallmatrix}\dots
$$
we come to the next M\"obius $p$-function $\tilde{\mu }_{5} \left(n\right)$ ($n$ runs only through the sifted numbers $(n,2,3)=1$ except 1)
\begin{equation*}
\begin{array}{l}
\tilde{\mu }_{5} \left(5\cdot n\right)=+1 \qquad \text{if 5$\nmid n$},\\ 
\tilde{\mu }_{5} \left(5\cdot n\right)=0 \qquad \;\;\, \text{if 5$\mid n$}.
\end{array}
\end{equation*}
The new values in the sequence above for prime 5 (5-sequence) are
$$
\setcounter{MaxMatrixCols}{30}
\begin{smallmatrix}
1&\quad\!\!&\quad\!\!&\quad\!\!&5&\quad\!\!&7&\quad\!\!&\quad\!\!&\quad\!\!&11&\quad\!\!&13&\quad\!\!&\quad\!\!&\quad\!\!&17&\quad\!\!&19&\quad\!\!&\quad\!\!&\quad\!\!&23&\quad\!\!&25&\quad\!\!&\quad\!\!&\quad\!\!&29&\\
+&&&&-&&-&&&& -&&- &&&&- &&- &&&&- &&0 &&&&- &
\end{smallmatrix}\dots
$$
 
Insert these new values into the Complete Sequence and making at the same time the following changes
\begin{equation*}
\begin{array}{l}
\tilde{\mu }_{5} \left(2\cdot 5\cdot n\right)=-1 \qquad\;\;\;\;\, \text{if 5$\nmid n$},\\ 
\tilde{\mu }_{5} \left(2\cdot 5\cdot n\right)=0 \qquad\;\;\;\;\;\;\,\, \text{if 5$\mid n$},\\
\tilde{\mu }_{5} \left(3\cdot 5\cdot n\right)=-1 \qquad\;\;\;\;\, \text{if 5$\nmid n$},\\
\tilde{\mu }_{5} \left(3\cdot 5\cdot n\right)=0 \qquad\;\;\;\;\;\;\,\, \text{if 5$\mid n$},\\
\tilde{\mu }_{5} \left(2\cdot 3\cdot 5\cdot n\right)=+1 \qquad \text{if 5$\nmid n$},\\
\tilde{\mu }_{5} \left(2\cdot 3\cdot 5\cdot n\right)=0 \qquad \;\;\, \text{if 5$\mid n$}.
\end{array}
\end{equation*}
For all other $n$ $\tilde{\mu }_{5} \left(n\right)$ remains unchanged
$$
\tilde{\mu }_{5} \left(n\right)=\tilde{\mu }_{3} \left(n\right).
$$
$$\setcounter{MaxMatrixCols}{30}
\begin{smallmatrix}
1&2&3&4&5&6&7&8&9&10&11&12&13&14&15&16&17&18&19&20&21&22&23&24&25&26&27&28&29&30\\
+&-&-&0&-&+&-&0&0& +& -& 0&- &+ &+ &0 &- &0 &- &0 &+ &+ &- &0 &0 &+ &0 &0 &- &-
\end{smallmatrix}
\dots
$$
All these transformations can be represented as
\small
\begin{equation} \label{Eq_19} 
\begin{array}{rcl}{\tilde{\mu }}_{5}{\left(n\right)} & {=} & {\tilde{\mu }_{3} \left(n\right)+2\left(1-\delta \left(\frac{n}{2} \right)\right)\left(1-\delta \left(\frac{n}{3} \right)\right)\delta \left(\frac{n-5^{2} }{5} \right)} \\ {} & {} & -2\left(1-\delta \left(\frac{n-2^{2} \cdot 5^{2} }{2^{2} } \right)\right)\left(1-\delta \left(\frac{n}{3} \right)\right)\delta \left(\frac{n-2\cdot 5^{2} }{2\cdot 5} \right)
\\ {} & {} & -2\left(1-\delta \left(\frac{n}{2} \right)\right)\left(1-\delta \left(\frac{n-3^{2} \cdot 5^{2} }{3^{2} } \right)\right)\delta \left(\frac{n-3\cdot 5^{2} }{3\cdot 5} \right)
\\ {} & {} & +2\left(1-\delta \left(\frac{n-2^{2} \cdot 5^{2} }{2^{2} } \right)\right)\left(1-\delta \left(\frac{n-2\cdot 3^{2} \cdot 5^{2} }{2\cdot 3^{2} } \right)\right)\delta \left(\frac{n-2\cdot 3\cdot 5^{2} }{2\cdot 3\cdot 5} \right)
\\ {} & {} & -\left(1-\delta \left(\frac{n}{2} \right)\right)\left(1-\delta \left(\frac{n}{3} \right)\right)\delta \left(\frac{n-5^{2} }{5^{2} } \right)
\\ {} & {} & +\left(1-\delta \left(\frac{n-2^{2} \cdot 5^{2} }{2^{2} } \right)\right)\left(1-\delta \left(\frac{n}{3} \right)\right)\delta \left(\frac{n-2\cdot 5^{2} }{2\cdot 5^{2} } \right)
\\ {} & {} & +\left(1-\delta \left(\frac{n}{2} \right)\right)\left(1-\delta \left(\frac{n-3^{2} \cdot 5^{2} }{3^{2} } \right)\right)\delta \left(\frac{n-3\cdot 5^{2} }{3\cdot 5^{2} } \right)
\\ {} & {} & -\left(1-\delta \left(\frac{n-2^{2} \cdot 5^{2} }{2^{2} } \right)\right)\left(1-\delta \left(\frac{n-2\cdot 3^{2} \cdot 5^{2} }{2\cdot 3^{2} } \right)\right)\delta \left(\frac{n-2\cdot 3\cdot 5^{2} }{2\cdot 3\cdot 5^{2} } \right). \end{array} 
\end{equation}
\normalsize 
For the asymptotic densities we have
\begin{equation} \label{Eq_20} 
\begin{array}{rcl} 
{\mathop{\lim }\limits_{n\to \infty } \frac{\tilde{M}_{5} \left(n\right)}{n} } & {=} & {\mathop{\lim }\limits_{n\to \infty } \left(\Pr \left(\tilde{\mu }_{5} \left(n\right)=+1\right)-\Pr \left(\tilde{\mu }_{5} \left(n\right)=-1\right)\right)} \\ {} & {=} & {-\left(1-\frac{1}{2} \right)^{2} \left(1-\frac{1}{3} \right)^{2} \left(1-\frac{1}{5} \right)^{2},} 
\end{array} 
\end{equation}
and 
\footnotesize
\begin{equation} \label{Eq_21} 
\mathop{\lim }\limits_{n\to \infty } \frac{\tilde{Q}_{5} \left(n\right)}{n} =\mathop{\lim }\limits_{n\to \infty } \left(\Pr \left(\left|\mu _{5} \left(n\right)\right|=1\right)\right)=\left(1-\frac{1}{2^{2} } \right)\left(1-\frac{1}{3^{2} } \right)\left(1-\frac{1}{5^{2} } \right).   
\end{equation}  
\normalsize

Denoting $p_{1} {\rm ,\; }p_{2} {\rm ,\; }p_{3} {\rm ,\; ...\; }$ the progression of primes $2,\, 3,\, 5\, ...$ we can continue above procedure until any $k^{{\rm th}} $ prime $p_{k} $ and get (by removing all multiples of $p_{1} {\rm ,\; }p_{2} {\rm ,\; }p_{3} {\rm ,\; ...\; ,\; }p_{k-1} $ from Complete Sequence) $p_{k} $-sequence and the corresponding M\"obius $p$-function $\tilde{\mu }_{p_{k} } \left(n\right)$ (below for each combination of primes $n$ runs only through the numbers of $p_{k} $-sequence $(n,p_{1} {\rm ,\; }p_{2} {\rm ,\; ...\; ,\; }p_{k-1} )=1$ except 1; all primes placed in ascending order with all possible distinct combinations, $m$ is position of primes, $1\le m\le k-1$)
\small
\begin{tabbing}
 \qquad\ \=$\tilde{\mu }_{p_{k} } \left(p_{k} \cdot n\right)=+1$ \qquad\qquad\qquad\qquad\qquad\quad\ \= if $p_{k}\nmid n$,\\ 
 \>$\tilde{\mu }_{p_{k} } \left(p_{k} \cdot n\right)=0$                                \> if $p_{k}\mid n$,\\    
 \>$\tilde{\mu }_{p_{k} } \left(p_{j_{k-1} } \cdot p_{k} \cdot n\right)=-1$            \> if $p_{k}\nmid n$, \;\! $1\le j_{k-1} \le k-1$,\\
 \>$\tilde{\mu }_{p_{k} } \left(p_{j_{k-1} } \cdot p_{k} \cdot n\right)=0$             \> if $p_{k}\mid n$, \;\! $1\le j_{k-1} \le k-1$,
\end{tabbing}
\noindent \qquad\ \dots
\small
\begin{equation} \label{Eq_22}
\tilde{\mu }_{p_{k} } \left(\underbrace{p_{j_{m} } \cdot ...\cdot p_{j_{k-1} } \cdot p_{k} }_{k-m+1={\rm number\; of\; primes}} \cdot n\right)=\left(-1\right)^{k-m+2} \;\, {\rm if} \; p_{k}\nmid n{\rm ,} \;\; 1\le j_{m} \le m{\rm ,}
\end{equation}
\begin{equation} \label{Eq_23}
\tilde{\mu }_{p_{k} } \left(p_{j_{m} } \cdot ...\cdot p_{j_{k-1} } \cdot p_{k} \cdot n\right)=0  \qquad\qquad\qquad\quad {\rm if} \; p_{k}\mid n{\rm ,} \;\; 1\le j_{m} \le m{\rm ,}
\end{equation}      
\noindent \qquad\ \dots
\begin{tabbing}
 \qquad\ \= $\tilde{\mu }_{p_{k} } \left(p_{j_{2} } \cdot ...\cdot p_{j_{k-1} } \cdot p_{k} \cdot n\right)=\left(-1\right)^{k}$ \qquad\qquad\;\; \= if $p_{k}\nmid n$, \;\! $1\le j_{2} \le 2$,\\
 \> $\tilde{\mu }_{p_{k} } \left(p_{j_{2} } \cdot ...\cdot p_{j_{k-1} } \cdot p_{k} \cdot n\right)=0$                          \> if $p_{k}\mid n$,     \;\! $1\le j_{2} \le  2$,\\
 \> $\tilde{\mu }_{p_{k} } \left(\underbrace{p_{1} \cdot p_{2} \cdot ...\cdot p_{k-1} \cdot p_{k} }_{{\rm all\; first\; }k{\rm \; primes}} \cdot       n\right)=\left(-1\right)^{k+1}$   \> if $p_{k}\nmid n$,\\
 \> $\tilde{\mu }_{p_{k} } \left(p_{1} \cdot p_{2} \cdot ...\cdot p_{k-1} \cdot p_{k} \cdot n\right)=0$  \> if $p_{k}\mid n$.
\end{tabbing}
\normalsize For all other $n$ $\tilde{\mu }_{p_{k} } \left(n\right)$ remains unchanged 
\[\tilde{\mu }_{p_{k} } \left(n\right)=\tilde{\mu }_{p_{k-1} } \left(n\right).\] 

Now compare what we obtained for $\tilde{\mu }_{p_{k} } \left(n\right)$ with the product
\footnotesize
\begin{equation} \label{Eq_24} 
\begin{array}{l} {\underbrace{-\left(1-\frac{1}{p_{1} } \right)^{2} ...\left(1-\frac{1}{p_{k-1} } \right)^{2} \left(1-\frac{1}{p_{k} } \right)^{2} }_{{\rm analog\; of\; }\tilde{\mu }_{p_{k} } \left(n\right)} =} \\ {\underbrace{-\left(1-\frac{1}{p_{1} } \right)^{2} ...\left(1-\frac{1}{p_{k-1} } \right)^{2} }_{{\rm analog\; of\; }\tilde{\mu }_{p_{k-1} } \left(n\right)} +\underbrace{\left(1-\frac{1}{p_{1} } \right)...\left(1-\frac{1}{p_{k-1} } \right)}_{{\rm analog\; of\; sieving\; factors}}} \\ 
{\times \left(\underbrace{2\left(1-\frac{1}{p_{1} } \right)...\left(1-\frac{1}{p_{k-1} } \right)\frac{1}{p_{k} } }_{{\rm analog\; of\; }\delta -{\rm terms\; with\; different\; combinations\; of\; distinct\; primes}} \right.} \\ 
{\left. -\underbrace{\left(1-\frac{1}{p_{1} } \right)...\left(1-\frac{1}{p_{k-1} } \right)\frac{1}{p_{k}^{2} } }_{{\rm analog\; of\; }\delta -{\rm terms\; reducing\; values\; of\; }\tilde{\mu }_{p_{k} } \left(n\right){\rm \; to\; zero}} \right).}
\end{array} 
\end{equation}
\normalsize
Put side by side the delta functions corresponding \eqref{Eq_22} and \eqref{Eq_23}
$$
\begin{array}{l} 
{\left(-1\right)^{k-m+2} 2\left(1-\delta \left(\frac{n-p_{1}^{2} }{p_{1} } \right)\right)...\left(1-\delta \left(\frac{n-p_{j_{m} }^{2} }{p_{j_{m} }^{2} } \right)\right)...} \\ {\times \left(1-\delta \left(\frac{n-p_{j_{k-1} }^{2} }{p_{j_{k-1} }^{2} } \right)\right)...\left(1-\delta \left(\frac{n-p_{k-1}^{2} }{p_{k-1} } \right)\right)\delta \left(\frac{n-p_{j_{m} } \cdot ...\cdot p_{j_{k-1} } \cdot p_{k}^{2} }{p_{j_{m} } \cdot ...\cdot p_{j_{k-1} } \cdot p_{k} } \right)} \\ {-\left(-1\right)^{k-m+2} \left(1-\delta \left(\frac{n-p_{1}^{2} }{p_{1} } \right)\right)...\left(1-\delta \left(\frac{n-p_{j_{m} }^{2} }{p_{j_{m} }^{2} } \right)\right)...} \\ {\times \left(1-\delta \left(\frac{n-p_{j_{k-1} }^{2} }{p_{j_{k-1} }^{2} } \right)\right)...\left(1-\delta \left(\frac{n-p_{k-1}^{2} }{p_{k-1} } \right)\right)\delta \left(\frac{n-p_{j_{m} } \cdot ...\cdot p_{j_{k-1} } \cdot p_{k}^{2} }{p_{j_{m} } \cdot ...\cdot p_{j_{k-1} } \cdot p_{k}^{2} } \right)} 
\end{array}
$$
and terms from the product
$$
\begin{array}{l} 
{\left(-1\right)^{k-m+2} \left(1-\frac{1}{p_{1} } \right)...\left(1-\frac{1}{p_{j_{m} } } \right)...\left(1-\frac{1}{p_{j_{k-1} } } \right)...\left(1-\frac{1}{p_{k-1} } \right)} \\ {\times \left(2\frac{1}{p_{j_{m} } \cdot ...\cdot p_{j_{k-1} } \cdot p_{k} } -\frac{1}{p_{j_{m} } \cdot ...\cdot p_{j_{k-1} } \cdot p_{k}^{2} } \right).} 
\end{array}
$$ 

We have a close analogy. For arbitrary term with the delta function we can find unique corresponding term in product \eqref{Eq_24} with the same law of parity. Reverse statement is also true. It means that we can get general expression for $\tilde\mu _{p_{k} } \left(n\right)$ through the delta functions very similar to the product \eqref{Eq_24}. Then applying the rules \eqref{Eq_16} and replacing the delta functions by corresponding frequencies or probabilities, we deduce the general formulae for the asymptotic densities
\footnotesize
\begin{equation} \label{Eq_25} 
\begin{array}{rcl} 
{\mathop{\lim }\limits_{n\to \infty } \frac{\tilde{M}_{p_{k} } \left(n\right)}{n} } & {=} & {\mathop{\lim }\limits_{n\to \infty } \left(\Pr \left(\tilde{\mu }_{p_{k} } \left(n\right)=+1\right)-\Pr \left(\tilde{\mu }_{p_{k} } \left(n\right)=-1\right)\right)} \\
{} & {=} & {-\left(1-\frac{1}{p_{1} } \right)^{2} \left(1-\frac{1}{p_{2} } \right)^{2} \cdot ...\cdot \left(1-\frac{1}{p_{k} } \right)^{2} } \\ 
{} & {=} & {-\prod _{m=1}^{k}\left(1-\frac{1}{p_{m} } \right)^{2},} 
\end{array} 
\end{equation} 
\begin{equation} \label{Eq_26} 
\begin{array}{rcl} 
{\mathop{\lim }\limits_{n\to \infty } \frac{\tilde{Q}_{p_{k} } \left(n\right)}{n} } & {=} & {\mathop{\lim }\limits_{n\to \infty } \left(\Pr \left(\left|\mu _{p_{k} } \left(n\right)\right|=1\right)\right)} \\ 
{} & {=} & {\left(1-\frac{1}{p_{1}^{2} } \right)\left(1-\frac{1}{p_{2}^{2} } \right)...\left(1-\frac{1}{p_{k}^{2} } \right)=\prod _{m=1}^{k}\left(1-\frac{1}{p_{m}^{2} } \right).} 
\end{array} 
\end{equation}
\normalsize 

From the way how we build the values of the M\"obius $p$-functions $\tilde{\mu }_{p_{k} } \left(n\right)$ and the definition \eqref{Eq_1} of $\mu \left(n\right)$ it is obvious that for $n<p_{k+1}^{2} $
\begin{equation} \label{Eq_27} 
\begin{array}{l} 
{\mu \left(n\right)=\tilde{\mu }_{p_{k} } \left(n\right),} \\ 
{M\left(n\right)=\tilde{M}_{p_{k} } \left(n\right).}
\end{array} 
\end{equation}
\\

\section{\textbf{Infinitude of primes}}

It is well known that there are an infinite number of primes (Euclid's theorem) and above procedure can be continued infinitely. However, it is important to know to what consequences for squarefree numbers we can come if we assume that there are a finite number of primes $p_{1} {\rm ,\; }p_{2} {\rm ,\; ...\; ,\; }p_{\max } $.

From that assumption and the definition \eqref{Eq_1} follows that the product of all primes $p_{1} \cdot p_{2} \cdot ...\cdot p_{\max } $ is the largest squarefree number - factorization of this number contains all available primes (greatest number); all other numbers are either squarefree with less quantity of primes in factorization or simply squareful. It means that the next bigger integer definitely should be squareful, i.e. $\mu \left(p_{1} \cdot p_{2} \cdot ...\cdot p_{\max } +1\right)=0$.  However, from the identity \eqref{Eq_3} follows that all delta functions for $\mu \left(p_{1} \cdot p_{2} \cdot ...\cdot p_{\max } +1\right)$ except the first one $-\delta \left(\frac{p_{1} \cdot p_{2} \cdot ...\cdot p_{\max } +1}{1\cdot 1} \right)=-1$ are 0 since the integer $p_{1} \cdot p_{2} \cdot ...\cdot p_{\max } +1$ and divisors in the delta functions do not have even a single common prime factor ($p_{1} \cdot p_{2} \cdot ...\cdot p_{\max } +1$ and  each divisor is coprime). Therefore $\mu \left(p_{1} \cdot p_{2} \cdot ...\cdot p_{\max } +1\right)=-1$, i.e. there is bigger squarefree number than $p_{1} \cdot p_{2} \cdot ...\cdot p_{\max } $. That is impossible - contradiction. Primes cannot be a finite number and there are infinitely many primes. This proof has some similarity with the Euclid's proof.

Let us consider an assumption that there are a finite number of primes  $p_{1} {\rm ,\; }p_{2} {\rm ,\; ...\; ,\; }p_{\max } $ from another point of view. We have seen already that in this case the largest squarefree number should be $p_{1} \cdot p_{2} \cdot ...\cdot p_{\max } $. After the largest squarefree number we observe the ``perfect calm'' - all integers bigger are squareful and, consequently, the asymptotic density of squarefree numbers is 0. However, this is not in agreement with formula \eqref{Eq_26}
\begin{equation} \label{Eq_28} 
\begin{array}{rcl} 
{\mathop{\lim }\limits_{n\to \infty } \frac{\tilde{Q}_{p_{\max } } \left(n\right)}{n} } & {=} & {\mathop{\lim }\limits_{n\to \infty } \left(\Pr \left(\left|\mu _{p_{\max } } \left(n\right)\right|=1\right)\right)} \\ {} & {=} & {\left(1-\frac{1}{p_{1}^{2} } \right)\left(1-\frac{1}{p_{2}^{2} } \right)...\left(1-\frac{1}{p_{\max }^{2} } \right)\ne 0.} 
\end{array} 
\end{equation} 

We can come to similar conclusion if we consider asymptotic density of prime numbers. It is clear that asymptotic density of finite number of primes is zero. For finite number of primes we have to cut infinite product \eqref{Eq_4} for prime detecting function on prime $p_{\max } $
\begin{equation} \label{Eq_29} 
\begin{array}{l} 
{p_{{\rm max}}^{\delta} \left(n\right)=\left(1-\delta \left(\frac{n-p_{1}^{2} }{p_{1} } \right)\right)\left(1-\delta \left(\frac{n-p_{2}^{2} }{p_{2} } \right)\right)} \\ {{\rm \; \; \; \; \; \; \; \; \; \; \; \; \; }...\left(1-\delta \left(\frac{n-p_{{\rm max}}^{2} }{p_{\max } } \right)\right), \qquad\qquad n\ge 2.}
\end{array} 
\end{equation} 
Then for cumulative sum of this function we get ``prime'' counting $p$-function
\begin{equation} \label{Eq_30} 
\pi _{p_{{\rm max}} } \left(n\right)=\sum _{m=2}^{n}p_{{\rm max}}^{\delta} \left(m\right).        
\end{equation}
Actually this function is counting primes and numbers which are relatively prime to the  primes $p_{1} {\rm ,\; }p_{2} {\rm ,\; ...\; ,\; }p_{\max } $ in the interval between 1 and $n$. Replacing the delta functions above by frequencies or probabilities we deduce that asymptotic densities of such integers are
\begin{equation} \label{Eq_31} 
\mathop{\lim }\limits_{n\to \infty } \frac{\pi _{p_{{\rm max}} } \left(n\right)}{n} =\left(1-\frac{1}{p_{1} } \right)\left(1-\frac{1}{p_{2} } \right)\cdot ...\cdot \left(1-\frac{1}{p_{{\rm max}} } \right)\ne 0. 
\end{equation} 
Again we are coming to contradiction. Consequently, our assumption is not appropriate and we have to conclude that the sequence of primes does not end.

\section{\textbf{Prime Number Theorem}}

The asymptotic law of the distribution of prime numbers is known as the Prime Number Theorem \cite{Hardy, Jameson, Tenenbaum, Weisstein2, Wikipedia2} -- if $\pi \left(n\right)$ is the number of primes $p\le n$, then $\pi \left(n\right)$ is  asymptotically equal to $\frac{n}{\ln n} $ as $n\to \infty $. For the first time this assertion was proved independently by Hadamard and de la Vallée Poussin~in 1896. In proofs they involved sophisticated complex analysis. Later the complex analytical proofs with various simplifications were found \cite{Jameson, Newman, Weisstein2, Wikipedia2}. An elementary proof was discovered by Selberg~and Erd\"os~in 1948. Their proof was without any use of complex analysis but it was not as clear as the analytic one. As noted by Hardy and Wright [\cite{Hardy}, p.~9], although it is elementary, "this proof is not easy". Despite the modern developments in the elementary methods \cite{Diamond, Wikipedia2}, proof of the Prime Number Theorem remains quite intricate.
 
We present our two versions of the proof, which is elementary and easy. Also we would emphasize that both proofs appears naturally without involving any ideas beyond the scope of traditional operations with integer numbers; it uses only some basic notions about the integers -- unique factorization into primes, sieving of primes and asymptotic densities for some type of integers.

\textbf{First proof}. We can cut the infinite product \eqref{Eq_4} for the prime detecting function on some $k^{{\rm th}} $ prime
\footnotesize
\begin{equation} \label{Eq_32} 
\begin{aligned} 
p_{{k}}^{\delta} \left(n\right) := \left(1-\delta \left(\frac{n-p_{1}^{2} }{p_{1} } \right)\right)\left(1-\delta \left(\frac{n-p_{2}^{2} }{p_{2} } \right)\right)\dots \left(1-\delta \left(\frac{n-p_{k}^{2} }{p_{k} } \right)\right)  
\end{aligned} 
\end{equation}
\normalsize
to get corresponding summatory function 
\begin{equation} \label{Eq_33} 
\pi _{p_{k} } \left(n\right):=\sum _{m=2}^{n}p_{{k}}^{\delta} \left(m\right).
\end{equation}
Further, we investigate asymptotic densities of these functions. Obviously, for 
$$\pi_{p_{1}} \left(n\right)=\sum _{m=2}^{n}p_{{1}}^{\delta} \left(m\right)=\sum _{\text{odd numbers $\leq n$}} 1$$
we have an asymptotic density of odd numbers
\begin{equation} \label{Eq_34}
\mathop{\lim }\limits_{n\to \infty } \frac{\pi _{p_{1} } \left(n\right)}{n}=\frac{1}{2}=\left(1-\frac{1}{p_{1} } \right) . 
\end{equation}
Using the following two remarkable properties of the delta-functions
\begin{equation} \label{Eq_35}
\begin{aligned}
 & \delta \left(\frac{n-p_{m_1}^{2} }{p_{m_1} } \right)\dots  \delta \left(\frac{n-p_{m_{i}}^{2} }{p_{m_{i}} } \right) \delta \left(\frac{n-p_{k}^{2} }{p_{k} } \right) = \\
 & \delta \left( {\frac{n}{p_{m_1} }} \right) \dots \delta \left( {\frac{n}{p_{m_{i}} }} \right) \delta \left(\frac{n-p_{k}^{2} }{p_{k} } \right) = \\ 
 & \delta \left( {\frac{n}{{{p_{m_1}} \dots {p_{m_{i}}} \cdot p_{k} }} - \left\lceil {\frac{{{p_k}}}{{{p_{m_1}} \dots {p_{m_{i}}}}}} \right\rceil } \right), \quad m_1< \ldots <m_i<k,
\end{aligned} 
\end{equation}
and
\begin{equation} \label{Eq_36}
\begin{aligned}
& \mathop {\lim}\limits_{n \to \infty } \frac{{\sum\limits_{m = 2}^n {\delta \left( {\frac{m}{p_{k} \cdot q} - r_k} \right)} }}{n} = \frac{1}{p_{k} \cdot q} = \\
& \frac{1}{p_{k}} \mathop {\lim}\limits_{n \to \infty } \frac{{\sum\limits_{m = 2}^n {\delta \left( {\frac{m}{q} - r_q} \right)} }}{n} \text{     (for any finite $r_q$) ,}
\end{aligned}
\end{equation}
it is not difficult to get recurrence relation (assuming that $\mathop{\lim }\limits_{n\to \infty } \frac{\pi _{p_{k-1} } \left(n\right)}{n}$ exists)
\begin{equation} \label{Eq_37}
\begin{aligned}
& \mathop{\lim }\limits_{n\to \infty } \frac{\pi _{p_{k} }  \left(n\right)}{n} =\\
 & \mathop{\lim }\limits_{n\to \infty } \frac { \sum\limits_{m = 2}^n {\left( {1 - \delta \left(\frac{m-p_{k}^{2} }{p_{k} } \right)} \right)\prod\limits_{l = 1}^{k-1} {\left( {1 - \delta \left( {\frac{m - p_{l}^{2}}{{{p_l}}} } \right)} \right)} } }{n} =\\
& \mathop{\lim }\limits_{n\to \infty } \frac{\pi _{p_{k-1} } \left(n\right)}{n}-\mathop{\lim }\limits_{n\to \infty } \frac { \sum\limits_{m = 2}^n {{ \delta \left(\frac{m-p_{k}^{2} }{p_{k} } \right)} \prod\limits_{l = 1}^{k-1} {\left( {1 - \delta \left( {\frac{m - p_{l}^{2}}{{{p_l}}} } \right)} \right)} } }{n} = \\
& \left( {1 - \frac{1}{{{p_k}}}} \right)\mathop{\lim }\limits_{n\to \infty } \frac{\pi _{p_{k-1} } \left(n\right)}{n} .
\end{aligned}
\end{equation}
From \eqref{Eq_34}, repeatedly applying the recurrence relation \eqref{Eq_37}, we obtain
\begin{equation} \label{Eq_38}
\mathop{\lim }\limits_{n\to \infty } \frac{\pi _{p_{k} } \left(n\right)}{n} =\left(1-\frac{1}{p_{1} } \right)\left(1-\frac{1}{p_{2} } \right)\cdot ...\cdot \left(1-\frac{1}{p_{k} } \right) ,
\end{equation}
which is valid for any $k \in \{1, 2, 3,\ldots\}$ .
In the limit of infinitely large $k$ in \eqref{Eq_38} from \eqref{Eq_4} and \eqref{Eq_5}, we conclude that asymptotic density of primes is
\begin{equation} \label{Eq_39}
\mathop{\lim }\limits_{n\to \infty } \frac{\pi _{p_\infty } \left(n\right)}{n} = \mathop{\lim }\limits_{n\to \infty } \frac{\pi \left(n\right)}{n} = \prod\limits_{\text{$p$ prime}} {\left( {1 - \frac{1}{p}} \right)} . 
\end{equation}
  
We can get similar result analysing asymptotic of some ratios of the Harmonic $H\left( n \right)$ and Harmonic-like series ${H_k}\left( n \right)$ (Harmonic series, where we are successively sifting (eliminating) reciprocals of all multiples of primes ${p_k},{\rm{   }} k = 1,{\rm{ }}2,{\rm{ }}...$). Define recursively for $k = 0, 1, 2, \ldots$ (as- \\sume also by default \small $H\left( 0 \right): = 0,\,{\rm{ }}{H_k}\left( 0 \right): = 0$ and ${H_k}\left( {\frac{n}{p}} \right): = {H_k}\left( {\left\lfloor {\frac{n}{p}} \right\rfloor } \right)$ ). \normalsize \\
\leftline{$ \qquad\qquad\qquad {H_0}\left( n \right)\,:\, = H\left( n \right) = \sum\limits_{m = 1}^n {\frac{1}{m}},$}\\
\leftline{$\qquad\qquad\qquad {H_1}\left( n \right)\,:\, = {H_0}\left( n \right) - \frac{1}{{{p_1}}}{H_0}\left( {\frac{n}{{{p_1}}}} \right) = H\left( n \right) - \frac{1}{{{p_1}}}{H_0}\left( {\frac{n}{{{p_1}}}} \right),$}\\
\leftline{$\begin{aligned} 
\qquad\qquad\qquad {H_2}\left( n \right): & = {H_1}\left( n \right) - \frac{1}{{{p_2}}}{H_1}\left( {\frac{n}{{{p_2}}}} \right) \\ 
&  = H\left( n \right) - \frac{1}{{{p_1}}}{H_0}\left( {\frac{n}{{{p_1}}}} \right) - \frac{1}{{{p_2}}}{H_1}\left( {\frac{n}{{{p_2}}}} \right), \\
\end{aligned}$}
\leftline{\qquad\qquad\qquad \dots}
\begin{equation} \label{Eq_40} 
\begin{aligned}
{H_k}\left( n \right): & = {H_{k - 1}}\left( n \right) - \frac{1}{{{p_k}}}{H_{k - 1}}\left( {\frac{n}{{{p_k}}}} \right) \qquad\quad \\  
& = H\left( n \right) - \sum\limits_{m = 1}^k  {\frac{1}{{{p_m}}}{H_{m - 1}}\left( {\frac{n}{{{p_m}}}} \right)}, \\ 
\end{aligned}
\end{equation}
\leftline{\qquad\qquad\qquad \dots}
\begin{equation} \label{Eq_41} 
\begin{aligned}
{H_\infty}\left( n \right): & = \mathop {\lim}\limits_{k \to \infty } \left({H_{k - 1}}\left( n \right) - \frac{1}{{{p_k}}}{H_{k - 1}}\left( {\frac{n}{{{p_k}}}} \right) \right)\\  
& = H\left( n \right) - \sum\limits_{m = 1}^\infty  {\frac{1}{{{p_m}}}{H_{m - 1}}\left( {\frac{n}{{{p_m}}}} \right)}=1 . \\ 
\end{aligned}
\end{equation}
In \eqref{Eq_41} we applied identity
\begin{equation} \label{Eq_42}
H\left( n \right) = 1 + \sum\limits_{k = 1}^\infty  {\frac{1}{{{p_k}}}{H_{k - 1}}\left( {\frac{n}{{{p_k}}}} \right)},
\end{equation} which is followed from the unique factorization into primes. For example\\
\footnotesize
$\begin{array}{l}
\begin{aligned}
H\left( {10} \right) & = 1 + \frac{1}{2}{H_0}\left( {\frac{{10}}{2}} \right) + \frac{1}{3}{H_1}\left( {\frac{{10}}{3}} \right) + \frac{1}{5}{H_2}\left( {\frac{{10}}{5}} \right) + \frac{1}{7}{H_3}\left( {\frac{{10}}{7}} \right) + \text {zero terms} \\ 
& = 1 + \frac{1}{2}\left( {1 + \frac{1}{2} + \frac{1}{3} + \frac{1}{4} + \frac{1}{5}} \right) + \frac{1}{3}\left( {1 + \frac{1}{3}} \right) + \frac{1}{5}\left( 1 \right) + \frac{1}{7}\left( 1 \right) .\\ 
\end{aligned}
\end{array}$
\normalsize
We can calculate asymptotic of the ratio of above functions to Harmonic series. We see that
\small
\begin{equation} \label{Eq_43}
\begin{aligned}
\mathop {\lim}\limits_{n \to \infty } \frac{{{H_1}\left( n \right)}}{{H\left( n \right)}} & = \mathop {\lim}\limits_{n \to \infty } \frac{{H\left( n \right)}}{{H\left( n \right)}} - \frac{1}{{{p_1}}}\mathop {\lim}\limits_{n \to \infty } \frac{{H\left( {\frac{n}{{{p_1}}}} \right)}}{{H\left( n \right)}} = 1 - \frac{1}{{{p_1}}}\mathop {\lim}\limits_{n \to \infty } \frac{{\log \frac{n}{{{p_1}}}}}{{\log n}} \\ 
& = 1 - \frac{1}{{{p_1}}}\mathop {\lim}\limits_{n \to \infty } \frac{{\log n - \log {p_1}}}{{\log n}} = 1 - \frac{1}{{{p_1}}} .\\ 
\end{aligned}
\end{equation}
\normalsize	
Looking on this result we can formulate a conjecture
\begin{equation} \label{Eq_44}
\mathop {\lim}\limits_{n \to \infty } \frac{{{H_k}\left( n \right)}}{{H\left( n \right)}} = \left( {1 - \frac{1}{{{p_1}}}} \right)\left( {1 - \frac{1}{{{p_2}}}} \right) \cdot ... \cdot \left( {1 - \frac{1}{{{p_k}}}} \right) .
\end{equation}
If assuming existence of $\mathop {\lim}\limits_{n \to \infty } \frac{{{H_{k - 1}}\left( n \right)}}{{H\left( n \right)}}$ we could get a recurrence relation
\begin{equation} \label{Eq_45}
\mathop {\lim}\limits_{n \to \infty } \frac{{{H_k}\left( n \right)}}{{H\left( n \right)}} = \left( {1 - \frac{1}{{{p_k}}}} \right)\mathop {\lim}\limits_{n \to \infty } \frac{{{H_{k - 1}}\left( n \right)}}{{H\left( n \right)}},
\end{equation}
then we shall come to the conclusion based on the result \eqref{Eq_43}, that asymptotic \eqref{Eq_44} holds for any natural numbers $k$.
From \eqref{Eq_40} it is not difficult to get that
\small
\begin{equation} \label{Eq_46}
\begin{aligned}
\mathop {\lim}\limits_{n \to \infty } \frac{{{H_k}\left( n \right)}}{{H\left( n \right)}} & = \mathop {\lim}\limits_{n \to \infty } \frac{{{H_{k - 1}}\left( n \right)}}{{H\left( n \right)}} - \frac{1}{{{p_k}}}\mathop {\lim}\limits_{n \to \infty } \frac{{{H_{k - 1}}\left( {\frac{n}{{{p_k}}}} \right)}}{{H\left( n \right)}} \\ 
&  = \mathop {\lim}\limits_{n \to \infty } \frac{{{H_{k - 1}}\left( n \right)}}{{H\left( n \right)}} - \frac{1}{{{p_k}}}\mathop {\lim}\limits_{n \to \infty } \frac{{{H_{k - 1}}\left( {\frac{n}{{{p_k}}}} \right)}}{{\log \frac{n}{{{p_k}}} + \log {p_k}}} \\ 
&  = \mathop {\lim}\limits_{n \to \infty } \frac{{{H_{k - 1}}\left( n \right)}}{{H\left( n \right)}} - \frac{1}{{{p_k}}}\mathop {\lim}\limits_{n \to \infty } \frac{{{H_{k - 1}}\left( {\frac{n}{{{p_k}}}} \right)}}{{H\left( {\frac{n}{{{p_k}}}} \right)\left( {1 + \frac{{\log {p_k}}}{{\log \frac{n}{{{p_k}}}}}} \right)}}\\
&  = \left( {1 - \frac{1}{{{p_k}}}} \right)\mathop {\lim}\limits_{n \to \infty } \frac{{{H_{k - 1}}\left( n \right)}}{{H\left( n \right)}} \,. \\ 
\end{aligned}
\end{equation}
\normalsize 
Thus, we have proved that \eqref{Eq_44} is valid for all positive integers $k$ and
\begin{equation} \label{Eq_47}
\mathop {\lim}\limits_{n \to \infty } \frac{{{H_\infty }\left( n \right)}}{{H\left( n \right)}} = \mathop {\lim}\limits_{n \to \infty } \frac{1}{{H\left( n \right)}} = \prod\limits_{\text{$p$ prime}} {\left( {1 - \frac{1}{p}} \right)}  = 0.
\end{equation}
Also we obtained that asymptotic of $\frac{{H\left( n \right)}}{{{H_k}\left( n \right)}}$ in the limit of infinitely large $k$ is Harmonic series $H\left( n \right)$.

Now consider the product $\frac{{{\pi _{{p_k}}}\left( n \right)}}{n} \cdot \frac{{H\left( n \right)}}{{{H_k}\left( n \right)}}$. Obviously, from \eqref{Eq_38} and \eqref{Eq_44} we have
\begin{equation} \label{Eq_48}
\mathop {\lim}\limits_{n \to \infty } \left( {\frac{{{\pi _{{p_k}}}\left( n \right)}}{n}\frac{{H\left( n \right)}}{{{H_k}\left( n \right)}}} \right) = 1{\rm{   }}{\rm{.}}
\end{equation}

Actually this limit is independent of $k$; it holds for any $k$. That is true also for infinitely large $k$. However in the limit of infinitely large $k$ in \eqref{Eq_48}, we have from \eqref{Eq_33}, \eqref{Eq_4} and \eqref{Eq_5} the prime counting function $\pi \left( n \right)$ instead of ${\pi _{{p_k}}}\left( n \right)$ and as we proved above \eqref{Eq_47} simply $H\left( n \right)$ instead of $\frac{{H\left( n \right)}}{{{H_k}\left( n \right)}}$
\begin{equation} \label{Eq_49}
\mathop {\lim}\limits_{n \to \infty } \left( {\frac{{{\pi _{{p_\infty }}}\left( n \right)}}{n}\frac{{H\left( n \right)}}{{{H_\infty }\left( n \right)}}} \right){\rm{ = }}\mathop {\lim}\limits_{n \to \infty } \left( {\frac{{\pi \left( n \right)}}{n}H\left( n \right)} \right){\rm{ = 1}}
\end{equation}
or $\pi \left(n\right)\sim \frac{n}{H\left(n\right)} \sim \frac{n}{\ln n}$ as $n\to \infty$.
The Prime Number Theorem is proved.

\textbf{Second proof}. In the framework of our approach, we also prove the assertion that the Mertens function $M\left(n\right)$ is of the lower order than $n$, i.e. $M\left(n\right)=o\left(n\right)$ as $n\to \infty $, which is equivalent to the asymptotic law of the distribution of prime numbers \cite{Diamond, Hardy, Jameson, Tenenbaum}.

Proceed in \eqref{Eq_25} and \eqref{Eq_26} to the limit of the infinitely large primes. From obtained result \eqref{Eq_47} we have 
\small
\begin{equation} \label{Eq_50} 
\begin{aligned}
{\mathop{\lim }\limits_{n\to \infty } \frac{\tilde{M}_{p_{\infty } } \left(n\right)}{n} } & =  {\mathop{\lim }\limits_{n\to \infty } \left(\Pr \left(\tilde{\mu }_{p_{\infty } } \left(n\right)=+1\right)-\Pr \left(\tilde{\mu }_{p_{\infty } } \left(n\right)=-1\right)\right)} \\ & = {-\mathop{\lim }\limits_{k\to \infty } \left(1-\frac{1}{p_{1} } \right)^{2} \left(1-\frac{1}{p_{2} } \right)^{2} \cdot ...\cdot \left(1-\frac{1}{p_{k} } \right)^{2} } \\ & =  {-\prod _{m=1}^{\infty }\left(1-\frac{1}{p_{m} } \right)^{2} =0,}
\end{aligned}
\end{equation}
\normalsize
\begin{equation} \label{Eq_51} 
\begin{aligned}
{\mathop{\lim }\limits_{n\to \infty } \frac{\tilde{Q}_{p_{\infty } } \left(n\right)}{n} } & = {\mathop{\lim }\limits_{n\to \infty } \left(\Pr \left(\left|\mu _{p_{\infty } } \left(n\right)\right|=1\right)\right)} \\ & = {\mathop{\lim }\limits_{k\to \infty } \left(1-\frac{1}{p_{1}^{2} } \right)\left(1-\frac{1}{p_{2}^{2} } \right)...\left(1-\frac{1}{p_{k}^{2} } \right)} \\ & = {\prod _{m=1}^{\infty }\left(1-\frac{1}{p_{m}^{2} } \right) =\frac{6}{\pi ^{2} } {\rm \; \; \; \; .}} 
\end{aligned} 
\end{equation}
We see that with advance to extremely large primes, the ``squarefree'' numbers with $\tilde{\mu }_{p} \left(n\right)=-1$ and $\tilde{\mu }_{p} \left(n\right)=+1$ occur with about equal frequencies. Same result will be with the M\"obius function $\mu \left( n \right)$ and actual squarefree numbers since for any positive integer $n$ yields
\begin{equation} \label{Eq_52}
\begin{aligned}
\mu \left(n\right) & = \mathop{\lim }\limits_{k\to \infty } \tilde{\mu }_{p_{k} } \left(n\right)=\tilde{\mu }_{p_{\infty} } \left(n\right), \\ M\left(n\right) & = \mathop{\lim }\limits_{k\to \infty } \tilde{M}_{p_{k} } \left(n\right)=\tilde{M}_{p_{\infty} } \left(n\right).
\end{aligned}
\end{equation}
Let us prove it. We arbitrarily choose some positive integer $n$ for which the M\"obius function has some value $\mu \left( n \right)$ and compare ${\tilde \mu _{{p_\infty }}}\left( n \right)$ and $\mu \left( n \right)$. Number $n$ has unique factorization (factors placed in ascending order) $n = p_1^{{\alpha _1}}p_2^{{\alpha _2}} \cdot ... \cdot p_i^{{\alpha _i}}$. From the way how we build the M\"obius $p$-functions (see section Sieve procedure), ${\tilde \mu _{{p_k}}}\left( n \right)$ becomes equal to $\mu \left( n \right)$ at ${p_k} \ge {p_i}$. Since ${p_i} < \left\{ {{p_{i + 1}},{\rm{ }}{p_{i + 2}},{\rm{ }}...,{\rm{ }}{p_\infty }} \right\}$ the M\"obius $p$-function ${\tilde \mu _{{p_\infty }}}\left( n \right)$ is certainly equal to $\mu \left( n \right)$ and since we did not impose any restrictions on the choice of the positive integer $n$, the equality is true for any other positive integers. For infinitely large primes the M\"obius $p$-function and Mertens $p$-function become identical to the original M\"obius function and Mertens function, respectively.  That means \eqref{Eq_52} is proved and then from \eqref{Eq_50} follows the Prime Number Theorem:
\begin{equation} \label{Eq_53} 
\mathop {\lim}\limits_{n \to \infty } \frac{{M\left( n \right)}}{n} = \mathop {\lim}\limits_{n \to \infty } \frac{{{{\tilde M}_{p_{\infty }}
}\left( n \right)}}{n} = 0
\end{equation}
or $M\left(n\right)=o\left(n\right)$ as $n\to \infty$ \cite{Hardy, Havil}.
Asymptotic densities of squarefree numbers with odd ($\mu \left(n\right)=-1$) and even ($\mu \left(n\right)=+1$) number of factors are same $\frac{3}{\pi ^{2} } $.

We confirm also the equivalence of $M\left( n \right) = o\left( n \right)$ and \\ $\sum\limits_{i = 1}^\infty  {\frac{{\mu \left( i \right)}}{i}}  = 0$. Recall formulas for the probabilities of the M\"obius function for arbitrary $n$ (see \cite{RAbrarov2} Eqs. (6) and (7) )
\begin{equation} \label{Eq_54}
\begin{array}{l}
 \Pr \left( {\mu \left( n \right) =  - 1} \right) = \frac{1}{2}\sum\limits_{i = 1}^{\sqrt n } {\frac{{\mu \left( i \right)}}{{{i^2}}}}  + \frac{1}{2}{\left( {\sum\limits_{i = 1}^{\sqrt n } {\frac{{\mu \left( i \right)}}{i}} } \right)^2}, \\ 
 \Pr \left( {\mu \left( n \right) =  + 1} \right) = \frac{1}{2}\sum\limits_{i = 1}^{\sqrt n } {\frac{{\mu \left( i \right)}}{{{i^2}}}}  - \frac{1}{2}{\left( {\sum\limits_{i = 1}^{\sqrt n } {\frac{{\mu \left( i \right)}}{i}} } \right)^2}. \\ 
\end{array}
\end{equation}
For asymptotic we have
\begin{equation} \label{Eq_55}
\begin{array}{l}
 \mathop {\lim}\limits_{n \to \infty } \frac{{{Q^ - }\left( n \right)}}{n} = \frac{1}{2}\sum\limits_{i = 1}^\infty  {\frac{{\mu \left( i \right)}}{{{i^2}}}}  + \frac{1}{2}{\left( {\sum\limits_{i = 1}^\infty  {\frac{{\mu \left( i \right)}}{i}} } \right)^2} = \frac{3}{{{\pi ^2}}} + \frac{1}{2}{\left( {\sum\limits_{i = 1}^\infty  {\frac{{\mu \left( i \right)}}{i}} } \right)^2}, \\ 
 \mathop {\lim}\limits_{n \to \infty } \frac{{{Q^ + }\left( n \right)}}{n} = \frac{1}{2}\sum\limits_{i = 1}^\infty  {\frac{{\mu \left( i \right)}}{{{i^2}}}}  - \frac{1}{2}{\left( {\sum\limits_{i = 1}^\infty  {\frac{{\mu \left( i \right)}}{i}} } \right)^2} = \frac{3}{{{\pi ^2}}} - \frac{1}{2}{\left( {\sum\limits_{i = 1}^\infty  {\frac{{\mu \left( i \right)}}{i}} } \right)^2}. \\ 
 \end{array}
\end{equation}
Equivalence of $M\left( n \right) = o\left( n \right)$ and $\sum\limits_{i = 1}^\infty  {\frac{{\mu \left( i \right)}}{i}}  = 0$ is followed from the identity
\begin{equation} \label{Eq_56}
\mathop {\lim}\limits_{n \to \infty } \frac{{M\left( n \right)}}{n} = \mathop {\lim}\limits_{n \to \infty } \frac{{{Q^ + }\left( n \right) - {Q^ - }\left( n \right)}}{n} =  - {\left( {\sum\limits_{i = 1}^\infty  {\frac{{\mu \left( i \right)}}{i}} } \right)^2}.
\end{equation}
\\


\begin{thebibliography}{14}
\bibitem{RAbrarov1}
R. M. Abrarov and S. M. Abrarov, \textit{Formulas for positive, negative and zero values of the M\"obius function}, 
\indent\verb+http://arxiv.org/PS_cache/arxiv/pdf/0905/0905.0294v1.pdf+
\smallskip
\bibitem{WolframDemo}
R. M. Abrarov and S. M. Abrarov, \textit{How the Superposition of the Periodic Pulsations of +1 and -1 Generates Values of the Mertens Function}, from \href{http://demonstrations.wolfram.com/HowTheSuperpositionOfThePeriodicPulsationsOf1And1GeneratesVa/}{\textcolor{blue}{\textit{The Wolfram Demonstrations Project}}}.
\smallskip
\bibitem{RAbrarov2}
R. M. Abrarov and S. M. Abrarov, \textit{Probabilistic interpretation of the M\"obius function identity and the Riemann Hypothesis},
\indent\verb+http://arxiv.org/PS_cache/arxiv/pdf/1002/1002.1682v1.pdf+
\smallskip
\bibitem{Diamond}
H. G. Diamond, \textit{``Elementary methods in the study of the distribution of prime numbers'', }Bull. Amer. Math. Soc. N. S. \textbf{7}, number 3 (1982), 553-589.
\smallskip 
\bibitem{Hardy}
G. H. Hardy and E. M. Wright, \textit{An introduction to the theory of numbers, 5th ed.}, Oxford University Press, Oxford, 1979.
\smallskip 
\bibitem{Havil}
J. Havil, \textit{Gamma: exploring Euler's constant}, Princeton University Press, 2003
\bibitem{Jameson}
G. J. O. Jameson, \textit{The Prime Number Theorem} , Cambridge University Press, 2003.
\smallskip 
\bibitem{Newman}
D. J. Newman, \textit{"Simple analytic proof of the prime number theorem"},~Amer. Math. Monthly~\textbf{87}: 693--696, (1980).
\smallskip 
\bibitem{Tenenbaum}
G. Tenenbaum and M. M. France, \textit{The Prime Numbers and Their Distribution}, AMS, 2001.
\smallskip  
\bibitem{Weisstein1}
Weisstein, Eric W. \textit{$"$M\"obius Function.$"$} From MathWorld - A Wolfram Web Resource. \\ \indent\verb+http://mathworld.wolfram.com/MoebiusFunction.html+
\smallskip
\bibitem{Weisstein2}
Weisstein, Eric W.~\textit{"Prime Number Theorem."} From~MathWorld - A Wolfram Web Resource. \\ \indent\verb+http://mathworld.wolfram.com/PrimeNumberTheorem.html+
\smallskip
\bibitem{Weisstein3}
Weisstein, Eric W. \textit{$"$Squarefree.$"$} From MathWorld - A Wolfram Web Resource. \\
\indent\verb+http://mathworld.wolfram.com/Squarefree.html+
\smallskip
\bibitem{Wikipedia1}
Wikipedia contributors,
\href{http://en.wikipedia.org/wiki/Mertens_function}{\textit{$"$Mertens function$"$,}} Wikipedia, The Free Encyclopedia.
\smallskip
\bibitem{Wikipedia2}
Wikipedia contributors,
\href{http://en.wikipedia.org/wiki/Prime_number_theorem}{\textit{$"$Prime Number Theorem$"$,}} Wikipedia, The Free Encyclopedia.
\\
\end{thebibliography}
\end{document}